\documentclass[journal]{IEEEtran}

\usepackage{amsmath,amsfonts,amssymb,amsopn,amsthm}
\newtheorem{theorem}{Theorem}
\newtheorem{corollary}[theorem]{Corollary}
\newtheorem{definition}[theorem]{Definition}
\newtheorem{proposition}[theorem]{Proposition}

\usepackage{dsfont,mathtools}
\usepackage{stmaryrd}
\SetSymbolFont{stmry}{bold}{U}{stmry}{m}{n}
\usepackage{altvars}
\usepackage{algorithmic}
\usepackage{algorithm}
\usepackage{array}
\usepackage[dvipsnames]{xcolor}
\usepackage{pgfplots}
\pgfplotsset{compat=1.17}
\usepackage[caption=false,font=normalsize,labelfont=sf,textfont=sf]{subfig}
\usepackage{textcomp}
\usepackage{relsize}
\usepackage{stfloats}
\usepackage{url}
\usepackage{verbatim}
\usepackage{graphicx}
\usepackage{multirow}
\graphicspath{{figs/}}
\usepackage{cite}
\usepackage[hidelinks]{hyperref}
\usepackage[capitalise]{cleveref}
\crefname{lemma}{Lemma}{Lemmas}
\Crefname{lemma}{Lemma}{Lemmas}
\crefname{equation}{}{}
\Crefname{equation}{Equation}{Equations}
\crefname{figure}{Fig.}{Figs.}
\usepackage{ifpdf}
\ifpdf
\usepackage{crossreftools}
\fi
\begin{document}

\title{Generalized Canonical Polyadic Tensor Decompositions\\with General Symmetry}

\author{%
    Alex Mulrooney,~\IEEEmembership{Student~Member,~IEEE,}
    and
    David Hong,~\IEEEmembership{Member,~IEEE}%
    \thanks{%
        A. Mulrooney and D. Hong are with the
        Department of Electrical and Computer Engineering,
        University of Delaware,
        Newark, DE, 19716 USA.
        Corresponding author: David Hong (email: hong@udel.edu).
    }%
}

\maketitle

\begin{abstract}
    Canonical Polyadic (CP) tensor decomposition is a workhorse algorithm for discovering underlying low-dimensional structure in tensor data.
    This is accomplished in conventional CP decomposition by fitting a low-rank tensor to data with respect to the least-squares loss.
    Generalized CP (GCP) decompositions generalize this approach by allowing general loss functions that can be more appropriate, e.g., to model binary and count data or to improve robustness to outliers.
    However, GCP decompositions do not explicitly account for any symmetry in the tensors, which commonly arises in modern applications.
    For example, a tensor formed by stacking the adjacency matrices of a dynamic graph over time will naturally exhibit symmetry along the two modes corresponding to the graph nodes.
    In this paper, we develop a symmetric GCP (SymGCP) decomposition that allows for general forms of symmetry, i.e., symmetry along any subset of the modes.
    SymGCP accounts for symmetry by enforcing the corresponding symmetry in the decomposition.
    We derive gradients for SymGCP that enable its efficient computation via all-at-once optimization with existing tensor kernels.
    The form of the gradients also leads to various stochastic approximations that enable us to develop stochastic SymGCP algorithms that can scale to large tensors.
    We demonstrate the utility of the proposed SymGCP algorithms with a variety of experiments on both synthetic and real data.
\end{abstract}


\section{Introduction}

\IEEEPARstart{T}{ensor} decomposition techniques are fundamental tools for the analysis of data tensors, i.e., for data organized as multi-dimensional arrays.
In particular, the canonical polyadic (CP) tensor decomposition \cite{hitchcock1927expression, carroll1970analysis, harshman1970foundations} provides a powerful tool for uncovering underlying low-rank signals in data.
CP decompositions do so by approximating the data with a sum of rank-$1$ \emph{component} tensors, yielding a Kruskal tensor.
In many cases, each component captures an underlying phenomenon in the data as a result of the uniquenss properties of these low-rank tensors.
CP decompositions have found applications in a wide range of fields, including signal processing, computer vision, neuroscience, and deep learning \cite{kolda2009tensor}.

Conventional CP decompositions fit the data with respect to the least-squares loss, i.e., they use squared error to measure the fit of the low-rank approximation to the data tensor.
A number of workhorse algorithms have been developed for this setting; see, e.g., \cite{kolda2009tensor} for an overview.
However, other (non-least-squares) loss functions can be more appropriate in some applications.
Of particular interest have been the nonnegative least-squares loss \cite{cichocki2009fast, shashua2005non, welling2001positive}, Poisson losses for count data \cite{chi2012tensors, hansen2015newton, ranadive2021all}, and Bernoulli losses for binary data \cite{miettinen2011boolean, wang2020learning}.
However, it is natural in various settings to use other losses beyond these choices, e.g., to better capture noise statistics or to obtain greater robustness to outliers.
As a result, several recent methods have emerged that allow more general loss functions to be used.
In particular, the Generalized CP (GCP) decomposition \cite{hong2020gcp}
fits a low-rank tensor to data
with respect to general user-specified loss functions.
An all-at-once gradient-based optimization approach was developed in \cite{hong2020gcp}.
They derive the gradients for GCP and observe that these gradients take the form of a sequence of Matricized Tensor Times Kahtri-Rao Products (MTTKRPS), for which there are highly optimized tensor kernels.
To scale the method to very large tensors, \cite{kolda2020sgf} developed stochastic gradients for GCP and employed a stochastic gradient method for optimization.
Other approaches include stochastic mirror descent algorithms based on a Bregman divergence suited for each chosen loss function \cite{pu2022stochastic, liu2025inertial}, and a second-order Gauss-Newton method \cite{vandecappelle2020second}.

Beyond the choice of loss function, another way to incorporate knowledge about the data tensor is to add constraints on the decomposition.
An important constraint, and the focus of our paper, is symmetry across the modes of the decomposition.
Symmetry naturally arises in modern data tensors, e.g., when one considers tensors corresponding to higher-order statistical moments.
Such tensors are symmetric across all their modes.
Another example of symmetry is a tensor formed by stacking the adjacency matrices of a dynamic graph over time.
Such tensors exhibit symmetry along the two modes corresponding to the graph nodes.
Naturally, one seeks a decomposition with matching symmetry.
Significant work has been done on this problem in the context of conventional CP decomposition with respect to the usual least-squares loss.
Kolda \cite{kolda2015nof} proposes a gradient-based optimization approach for symmetric CP with respect to least-squares losses (including nonnegative least-squares).
More recently, \cite{sherman2020estimating} proposed a specialized algorithm for symmetric CP from moment tensors that exploits the tensor structure to obtain highly efficient implicit computation of the function values and gradients.
Another recent approach is the subspace power method \cite{kileel2019subspace}, which sequentially finds the best rank-$1$ tensor in subspaces of the matricized data tensor.
Thus far, an overall focus has been on computing (least-squares) symmetric CP decompositions, where the symmetry is across all the modes of the tensor.

In this paper, we develop algorithms for CP decompositions that allow for both general loss functions (beyond the usual least-squares) and general forms of symmetry (beyond the usual symmetry across all the tensor modes).
For this purpose, we first define a general notion of tensor symmetry that encompasses tensors whose entries are equal under permutations of subsets of their indices.
This generalizes the usual notion of tensor symmetry, where entries are equal under permutations of all the indices.
We then formulate our proposed Symmetric GCP (SymGCP) method and develop an all-at-once gradient-based optimization approach for fitting the decomposition to data.
To scale the method to large tensors, we finally develop an efficient stochastic gradient method analogous to \cite{kolda2020sgf}.
We demonstrate the utility of the proposed SymGCP algorithms on synthetic and real datasets that encompass a variety of data types including binary and count data.

\Cref{sec:background} reviews relevant notation and background for tensors.
\Cref{sec:general-symmetry} defines the general form of tensor symmetry we consider in this paper.
\Cref{sec:formulation} formulates the proposed SymGCP decomposition.
\Cref{sec:gradients} then derives formulas for the corresponding gradients that enable efficient computation,
and \cref{sec:stochastic-gradients} derives stochastic gradients that enable further scalability via stochastic optimization.
Finally, \cref{sec:numerical-experiments,sec:real-data-experiments} illustrate SymGCP through experiments on synthetic and real data,
and \cref{sec:conclusion} provides concluding remarks.

\section{Tensor Notations} \label{sec:background}
This section establishes notations and operations for tensors and tensor decompositions, which we will use throughout the remainder of the paper.

We use lowercase unbolded letters (e.g., $z$) to denote scalars,
lowercase bold letters (e.g., $\bmy$) to denote vectors,
uppercase bold letters (e.g., $\bmA$) to denote matrices,
and calligraphic bold letters (e.g., $\bclX$) to denote tensors.
We write $\bclX \in \bbR^{I_1 \times \dots \times I_N}$ to signify that $\bclX$ is a real-valued tensor with $N$ modes, whose size along mode $n$ is $I_n$.
We denote the entry of a tensor $\bclX$ at index $(i_1, \dots, i_N)$ as $x_{i_1, \dots, i_N}$.
We also sometimes use the multi-index $i = (i_1,\dots,i_N)$ to denote entries as $x_i$.

We denote the vectorization of a tensor $\bclX \in \bbR^{I_1 \times \dots \times I_N}$ as
\begin{equation}
    \operatorname{vec}(\bclX) =
    \begin{bmatrix}
        x_{1,1,\dots,1} \\
        x_{2,1,\dots,1} \\
        \vdots \\
        x_{I_1,1,\dots,1} \\
        x_{1,2,\dots,1} \\
        \vdots \\
        x_{I_1,I_2,\dots,I_N}
    \end{bmatrix} \in \bbR^{I_1 \cdots I_N}
\end{equation}
The fibers of a tensor $\bclX$ along a mode $n$ generalize the notion of rows and columns to higher-order tensors; each fiber is a vector obtained by collecting the entries along the given mode, i.e., they are slices of the form $\bclX_{i_1,\dots,i_{n-1},:,i_{n+1},\dots,i_N} \in \bbR^{I_n}$.
The mode-$n$ matricization of a tensor $\bclX$, denoted by $\bmX_{(n)}$, is a matrix of size $I_n \times \prod_{j \neq n}I_j$ formed by horizontally stacking all the mode-$n$ fibers of $\bclX$ as columns.

We use $\circ$ to denote the outer product,
$\langle \cdot, \cdot \rangle$ to denote the inner product,
$\otimes$ to denote the Kronecker product,
$\odot$ to denote the Khatri-Rao (i.e., column-wise Kronecker) product,
and $\ast$ to denote the Hadamard (i.e., entrywise) product.
We use $[n]$ to denote the set $\{1,2,\dots,n\}$ for positive integers $n$.

Finally, we denote Kruskal tensors as
\begin{equation}
    \label{eq:ktensor}
    \llbracket \bmlambda ; \bmA_{1}, \dots, \bmA_{N} \rrbracket
    =
    \sum_{j=1}^r \lambda_j \cdot ( \bmA_1 )_{:,j} \circ \cdots \circ ( \bmA_N)_{:,j}
    ,
\end{equation}
where $\bmlambda \in \bbR^r$ is the weight vector,
and $\bmA_n \in \bbR^{I_n \times r}$ is the mode-$n$ factor matrix for each $n \in [N]$.
This low-rank tensor is the output of CP decomposition.

\section{Symmetric Tensors with General Symmetry} \label{sec:general-symmetry}

This section describes the general form of tensor symmetry that we consider to set the stage for our proposed symmetric GCP method (\cref{sec:formulation}).
\Cref{sec:symmetry:full} describes the full tensor symmetry that is commonly used,
\cref{sec:symmetry:general} describes the more general form of tensor symmetry we consider,
and \cref{sec:symmetric-kruskal-tensors} describes the corresponding notion of symmetric Kruskal tensors.

\subsection{Full Tensor Symmetry}
\label{sec:symmetry:full}

Before we present the general form of tensor symmetry we consider, here we briefly review the conventional form of tensor symmetry that is commonly studied.
In particular, the term ``symmetric tensor'' often refers to a tensor that is equivalent under any permutation of \emph{all its indices}.
For an $N$-way tensor $\bclX \in \bbR^{I_1 \times \cdots \times I_N}$, this means that
\begin{equation}
    \forall_{\pi \in \Pi([N])} \;
    \forall_{i \in [I_1] \times \cdots \times [I_N]}
    \quad
    x_{\pi(i)} = x_i
    ,
\end{equation}
where $\Pi([N])$ denotes the symmetric group on $[N]$, i.e., the set of all permutations of $[N]$, and the permutation $\pi \in \Pi([N])$ permutes the tensor indices $i = (i_1, \dots, i_N)$ as
\begin{equation}
    \pi(i) = (i_{\pi_1}, \dots, i_{\pi_N})
    .
\end{equation}
This can be written equivalently as
\begin{equation}
    \forall_{\pi \in \Pi([N])} \quad
    \operatorname{permutedims}(\bclX, \pi) = \bclX
    ,
\end{equation}
where $\operatorname{permutedims}$ permutes the modes of $\bclX$ by $\pi$.
In the language of group theory, this mode-wise permutation is the group action of $\pi$ on $\bclX$ that is relevant for our context, and symmetry here means that $\bclX$ is invariant under $\Pi([N])$.

Note that this notion of symmetry requires that $I_1 = \cdots = I_N$. In this paper we will refer to such tensors as being ``fully symmetric'' to distinguish from the more general form of symmetry we describe next.

\subsection{General Tensor Symmetry}
\label{sec:symmetry:general}

We now define the general notion of tensor symmetry that will be our focus.
In particular, motivated by applications such as dynamic graphs, we seek a notion of tensor symmetry that allows for symmetry in just subsets of the modes.
Note that a tensor formed by stacking symmetric adjacency matrices along a third mode would be symmetric across the first two modes but not across all three.
Thus, we relax the requirement that the entries of the tensor be equal under any permutation of \emph{all its indices} to instead just be equal under any permutation of some specified \emph{subsets of its indices}.

As a concrete example, suppose we have a three-way tensor $\bclX \in \bbR^{I_1 \times I_2 \times I_3}$ that is symmetric in its first two modes.
This means that
\begin{equation}
    \forall_{i \in [I_1] \times [I_2] \times [I_3]}
    \quad
    x_{i_1, i_2, i_3} = x_{i_2, i_1, i_3}
    .
\end{equation}
In other words, the frontal slices $\bclX_{:,:,i_3}$ for $i_3 \in [I_3]$ are all symmetric matrices.
This symmetry can be systematically encoded by enumerating all the permutations under which the tensor is unchanged,
which in this case yields the permutations $\{(1,2,3), (2,1,3)\}$.
Namely, we allow any permutation of modes $1$ and $2$, but fix mode $3$.

Put another way, we partition the modes here into subsets of ``equivalent'' modes that may be arbitrarily permuted without changing the tensor.
For our example, we have the partition $\{1,2\} \sqcup \{3\}$ of the set $\{1,2,3\}$ of all the modes since the tensor is symmetric across modes $1$ and $2$.
The set of all permutations of modes $1$ and $2$ are now given by the symmetric group $\Pi(\{1,2\})$ and the (trivial) set of all permutations of mode $3$ are given by the symmetric group $\Pi(\{3\})$.
Combining them then yields the group of relevant permutations
\begin{equation}
    \Pi(\{1,2\}) \times \Pi(\{3\}) = \{(1,2,3), (2,1,3)\}
    ,
\end{equation}
where $\times$ denotes the direct product.
In the language of group theory, this is a permutation group and a subgroup of $\Pi([N])$.

Abstracting our example leads naturally to the following general notion of tensor symmetry.

\begin{definition}[Tensor Symmetry]
    \label{def:symmetry}
    A tensor $\bclX \in \bbR^{I_1 \times \cdots \times I_N}$ is symmetric with respect to a partition $\clI_1 \sqcup \cdots \sqcup \clI_K = [N]$ of its modes into $K$ cells if
    \begin{equation}
        \forall_{\pi \in \Pi(\clI_1) \times \cdots \times \Pi(\clI_K)}
        \quad
        \forall_{i \in [I_1] \times \cdots \times [I_N]}
        \quad
        x_{\pi(i)} = x_i
        ,
    \end{equation}
    where $\Pi(\clI_1) \times \cdots \times \Pi(\clI_K) \subseteq \Pi([N])$ is the group of all permutations that permute only indices within each cell of the partition.
    Stated in terms of mode-wise permutation, we have
    \begin{equation}
        \forall_{\pi \in \Pi(\clI_1) \times \cdots \times \Pi(\clI_K)} \quad
        \operatorname{permutedims}(\bclX, \pi) = \bclX
        .
    \end{equation}
    This definition requires the dimensions in each cell to match, i.e., $I_j = I_k$ for any $j,k \in \clI_{\ell}$ and $\ell \in [K]$.
\end{definition}

This general notion of tensor symmetry encompasses several previously studied cases:
\begin{itemize}
    \item \textbf{Fully symmetric} tensors correspond to a partition with one cell containing all of the indices, i.e., $\clI_1 = [N]$.
    \item \textbf{Nonsymmetric} tensors correspond to a partition with one cell for each index, i.e., $\clI_j = \{j\}$ for all $j \in [N]$.
    \item \textbf{INDSCAL} \cite{carroll1970analysis} enforces symmetry along the first two modes of a three-way tensor and corresponds to the partition $\clI_1 \sqcup \clI_2 = \{ 1,2 \} \sqcup \{3\}$.
\end{itemize}

\subsection{Symmetric Kruskal Tensors with General Symmetry} \label{sec:symmetric-kruskal-tensors}

We now define a symmetric analogue of the (nonsymmetric) Kruskal tensor in \cref{eq:ktensor} for the general notion of tensor symmetry (Definition~\ref{def:symmetry}) that is our focus.
The overall idea is to constrain the Kruskal tensor to have the appropriate symmetry by simply constraining the factor matrices in each cell to match.
Namely,
given a partition $\clI_1 \sqcup \cdots \sqcup \clI_K = [N]$ of the modes, we consider Kruskal tensors of the form
\begin{equation}
    \llbracket \bmlambda ; \bmA_{1}, \dots, \bmA_{N} \rrbracket
    \quad
    \text{s.t.}
    \quad
    \forall_{\ell \in [K]} \;
    \forall_{j,k \in \clI_\ell} \;\;
    \bmA_j = \bmA_k
    .
\end{equation}
Effectively, rather than having separate factor matrices for each mode, we have just a single factor matrix for each cell in the partition.
This leads naturally to the following notion of symmetric Kruskal tensors with a corresponding notation.

\begin{definition}[Symmetric Kruskal Tensors]
    \label{def:symktensor}
    A symmetric Kruskal tensor with respect to a partition $\clI_1 \sqcup \cdots \sqcup \clI_K = [N]$ of $N$ modes is given as
    \begin{align}
        &
        \llbracket \bmlambda ; \clI_1 \to \bmA_{1}, \dots, \clI_K \to \bmA_{K} \rrbracket
        =
        \llbracket \bmlambda ; \bmA_{\sigma_1}, \dots, \bmA_{\sigma_N}\rrbracket
        \\&\nonumber
        \qquad\qquad
        =
        \sum_{j=1}^r
        \lambda_j \cdot (\bmA_{\sigma_1})_{:,j} \circ \cdots \circ (\bmA_{\sigma_N})_{:,j}
        ,
    \end{align}
    where $\bmlambda \in \bbR^r$ is the weight vector,
    each $\bmA_k \in \bbR^{I_k \times r}$ is the factor matrix associated to the cell $\clI_k$,
    and each $\sigma_n$ identifies which cell contains mode $n$,
    i.e.,
    \begin{equation}
        \forall_{k \in [K]} \;
        \forall_{n \in \clI_k}
        \quad
        \sigma_n = k
        .
    \end{equation}
\end{definition}

Note that by construction, any symmetric Kruskal tensor
$\llbracket \bmlambda ; \clI_1 \to \bmA_{1}, \dots, \clI_K \to \bmA_{K}\rrbracket$
is symmetric with respect to the partition $\clI_1 \sqcup \cdots \sqcup \clI_K$.
The following proposition makes this statement precise.
It follows straightforwardly by substitution and simple rewriting;
we provide an elementary proof in \cref{sec:symmetric:ktensor:issymmetric}
for the reader's convenience.

\begin{proposition}[Symmetric Kruskal Tensors are Symmetric]
    \label{thm:symmetric:ktensor:issymmetric}
    For any partition $\clI_1 \sqcup \cdots \sqcup \clI_K$ of $N$ modes, factor matrices $\bmA_1 \in \bbR^{I_1 \times r}, \dots, \bmA_K \in \bbR^{I_K \times r}$, and weight vector $\bmlambda \in \bbR^r$, the symmetric Kruskal tensor $\llbracket \bmlambda ; \clI_1 \to \bmA_{1}, \dots, \clI_K \to \bmA_{K} \rrbracket$ is symmetric with respect to $\clI_1 \sqcup \cdots \sqcup \clI_K$.
\end{proposition}

\section{Problem Statement: Symmetric GCP (SymGCP)}
\label{sec:formulation}

We are now ready to state our proposed Symmetric GCP (SymGCP) decomposition problem.
The overall idea is to fit low-rank model tensors to data with respect to general losses while also constraining the model to exhibit the desired symmetry.
In particular, we constrain the model tensor to be a symmetric Kruskal tensor as described in Definition~\ref{def:symktensor}.
This leads to the following SymGCP problem formulation
\begin{equation}
    \label{eq:symgcp:problem}
    \min_{\bmlambda, \bmA_1, \dots, \bmA_K}
    \;
    \underbrace{
        \clL(
            \bclX,
            \llbracket
            \bmlambda ; \clI_1 \to \bmA_{1}, \dots, \clI_K \to \bmA_{K}
            \rrbracket
        )
    }_{= \clF(\bmlambda, \bmA_1, \dots, \bmA_K)}
    ,
\end{equation}
where
\begin{itemize}
    \item $\bclX \in \bbR^{I_1 \times \cdots \times I_N}$ is an $N$-way data tensor,
    \item $\clI_1 \sqcup \cdots \sqcup \clI_K = [N]$ partitions the $N$ modes,
    \item $\bmA_1 \in \bbR^{\tlI_1 \times r},\dots,\bmA_K \in \bbR^{\tlI_K \times r}$ are the corresponding factor matrices with sizes defined for each $k \in [K]$ as
        \begin{equation*}
            \tlI_k = \text{the unique element of } \{I_n : n \in \clI_k\}
            ,
        \end{equation*}
    \item $\bmlambda \in \bbR^r$ are the desired weights,
        and
    \item $\clL : \bbR^{I_1 \times \cdots \times I_N} \times \bbR^{I_1 \times \cdots \times I_N} \to \bbR$ is a general loss function defined between the data tensor $\bclX$ and the model tensor
        $\bclM = \llbracket
        \bmlambda ; \clI_1 \to \bmA_{1}, \dots, \clI_K \to \bmA_{K}
        \rrbracket$.
\end{itemize}
Note that the dimensions of the data tensor $\bclX$ must respect the symmetry defined by the partition $\clI_1 \sqcup \cdots \sqcup \clI_K$ for SymGCP to make sense, i.e., the set $\{I_n : n \in \clI_k\}$ must have only one unique element for each $k \in [K]$.

For the loss function $\clL$, we will focus on loss functions defined in terms of an entrywise loss $\ell : \bbR \times \bbR \to \bbR$ as
\begin{equation}
    \label{eq:loss:entrywise}
    \clL(\bclX, \bclM)
    =
    \sum_{i} \ell(x_i, m_i)
    ,
\end{equation}
where the sum is usually taken over all the entries of the tensor.
Entrywise weights can be incorporated straightforwardly as
\begin{equation}
    \clL(\bclX, \bclM)
    =
    \sum_{i} w_i \ell(x_i, m_i)
    ,
\end{equation}
where $\bclW \in \bbR^{I_1 \times \cdots \times I_N}$ is an associated weight tensor.
Weighted formulations make it possible to easily handle missing data by simply giving unobserved entries a weight of zero.
In the context of symmetric tensors, weights can also be used to construct losses that avoid ``double-counting'' the duplicated entries arising from symmetry; see, e.g., \cite[Section~4.2]{kolda2015nof} for more discussion.

We conclude this section by noting that the objective in \cref{eq:symgcp:problem} has a scaling ambiguity inherent to the Kruskal tensor, since
\begin{align*}
    &
    \llbracket
    \bmlambda ; \clI_1 \to \bmA_{1}, \dots, \clI_K \to \bmA_{K}
    \rrbracket
    \\&\qquad
    =
    \llbracket
    \rho^N \bmlambda ; \clI_1 \to \bmA_{1}/\rho, \dots, \clI_K \to \bmA_{K}/\rho
    \rrbracket
    ,
\end{align*}
for any $\rho > 0$.
As discussed in \cite[Section~4.4]{kolda2015nof},
the same ambiguity also appears in the context of the least-squares loss
(indeed, it appears for any loss).
Moreover, we found that this ambiguity caused issues in practice for optimization.
Here, we address this ambiguity similarly to \cite{kolda2015nof} through the addition of the following regularizer:
\begin{equation}
    r_{\gamma}(\bmA_1, \dots, \bmA_K)
    =
    \gamma
    \cdot
    \sum_{k=1}^K \sum_{j=1}^r ( \| (\bmA_k)_{:,j} \|^2 - 1 )^2,
\end{equation}
where $\gamma > 0$ is the regularization parameter.
This regularizer eliminates the scaling ambiguity by penalizing factor matrices whose columns are not unit norm.

\section{Gradients for Symmetric GCP} \label{sec:gradients}

This section derives gradients for the SymGCP problem \cref{eq:symgcp:problem} that show how it can be computed using existing tensor kernels, enabling efficient SymGCP decomposition via gradient-based optimization.
For simplicity, we focus here on loss functions $\clL$ defined entrywise as in \cref{eq:loss:entrywise};
incorporating weights or accommodating general loss functions can be done in a similar way.

The following theorem provides simple formulas for the SymGCP gradients.
See \cref{theorem:symgcp-gradients:proof} for its proof.

\begin{theorem}[SymGCP Gradients] \label{theorem:symgcp-gradients}
    Let $\clL(\bclX, \bclM)$ be a loss function defined by an entrywise loss $\ell(x,m)$ as in \cref{eq:loss:entrywise}.
    Then the SymGCP objective $\clF(\bmlambda, \bmA_1, \dots, \bmA_K)$ in \cref{eq:symgcp:problem} has gradients given by
    \begin{align}
        \label{eq:symgcp:grad:lam}
        \frac{\partial \clF}{\partial \bmlambda}
        &=
        (\odot_{n=1}^N \bmA_{\sigma_n})^T \operatorname{vec}\left(\bclY\right)
        , \\
        \label{eq:symgcp:grad:Ak}
        \frac{\partial \clF}{\partial \bmA_k}
        &=
        \sum_{t \in \clI_k} \bmY_{(t)} (\odot_{n \neq t} \bmA_{\sigma_n}) \operatorname{diag}(\bmlambda)
        ,
    \end{align}
    where
    \begin{align*}
        \odot_{j=1}^N \bmA_{\sigma_j}
        &= \bmA_{\sigma_N} \odot \cdots \odot \bmA_{\sigma_1}
        , \\
        \odot_{j \neq t} \bmA_{\sigma_j}
        &= \bmA_{\sigma_N} \odot \cdots \odot \bmA_{\sigma_{t+1}} \odot \bmA_{\sigma_{t-1}} \odot \cdots \odot \bmA_{\sigma_1}
        ,
    \end{align*}
    and $\bclY$ is the derivative tensor whose entries are
    \begin{equation*}
        y_i = \frac{\partial \ell}{\partial m_i} (x_i, m_i)
        .
    \end{equation*}
\end{theorem}

This theorem generalizes the gradient formulas for (least-squares) symmetric CP given in \cite[Section 4.3]{kolda2015nof} to general losses.
Likewise, this theorem generalizes the gradient formulas for nonsymmetric GCP given in \cite[Theorem~3]{hong2020gcp}.
Similar to (nonsymmetric) GCP, this theorem reveals that the gradients can be efficiently computed through a sequence of Matricized Tensor Times Khatri-Rao Products (MTTKRPs) followed by a simple aggregation step across the modes in each cell.
This enables efficient computation of the gradient through the use of optimized tensor kernels for these operations; see, e.g., \cite{phan2013fast}.

Note that \cref{theorem:symgcp-gradients} holds regardless of the symmetry of $\bclX$, i.e., it need not have any symmetry.
Consequently, these formulas can be used in settings where we expect symmetry in the latent phenomenon, but have a data tensor that is not symmetric due, e.g., to noise.

When the data tensor $\bclX$ is in fact symmetric with respect to the partition $\bclI_1 \sqcup \cdots \sqcup \bclI_K$ (as is indeed commonly the case), a further simplification of the gradients is available.
It relies on the following property of symmetric MTTKRPs.

\begin{proposition}[Symmetric MTTKRPs] \label{lem:symmetric-mttkrps}
    If $\bclY \in \bbR^{I_1 \times \cdots \times I_N}$ is symmetric with respect to the partition $\bclI_1 \sqcup \cdots \sqcup \bclI_K$,
    then for every $k \in [K]$, we have
    \begin{equation*}
        \forall_{u,v \in \clI_k}
        \quad
        \bmY_{(u)}(\odot_{j \neq u} \bmA_{\sigma_j})
        =
        \bmY_{(v)}(\odot_{j \neq v} \bmA_{\sigma_j})
        ,
    \end{equation*}
    where $\sigma$ is as in Definition~\ref{def:symktensor}.
\end{proposition}

We provide a simple proof in \cref{lem:symmetric-mttkrps:proof}.
Applying this proposition to the gradient formula \cref{eq:symgcp:grad:Ak} immediately yields the following simplified form of the factor matrix gradients for the case where $\bclX$ is symmetric.

\begin{corollary}[SymGCP Gradients for Symmetric Data]
    If $\bclX \in \bbR^{I_1 \times \cdots \times I_N}$
    is symmetric with respect to the partition $\bclI_1 \sqcup \cdots \sqcup \bclI_K$,
    then the gradient \cref{eq:symgcp:grad:Ak} of the SymGCP objective $\clF(\bmlambda, \bmA_1, \dots, \bmA_K)$ with respect to factor matrix $\bmA_k$ simplifies as
    \begin{equation}
        \frac{\partial \clF}{\partial \bmA_k}
        =
        | \clI_k | \bmY_{(k')} (\odot_{j \neq k'} \bmA_{\sigma_j}) \operatorname{diag} (\bmlambda)
        ,
    \end{equation}
    where $| \clI_k |$ denotes the size of $\clI_k$,
    and $k' \in \clI_k$ is an arbitrary mode from $\clI_k$.
\end{corollary}

This form enables even more efficient computation of the SymGCP gradients since it reduces the number of MTTKRPs needed from $N$ to $K$. It also replaces the summation of the (nonsymmetric) gradients with a simple scaling.

\section{Stochastic gradients for Symmetric GCP} \label{sec:stochastic-gradients}

For very large tensors,
computing the full gradients derived in \cref{sec:gradients} can be prohibitively expensive.
Here, we develop stochastic approximations that can enable scalable SymGCP via stochastic optimization.

In particular, we consider stochastic gradients analogous to those proposed in \cite{kolda2020sgf} for nonsymmetric GCP.
The approach is motivated by the observation that for a data tensor $\bclX$ with $n^N$ entries, computing the full gradient with respect to any of the factor matrices or the weights requires computing and storing the derivative tensor $\bclY$, which is the same size as $\bclX$.
The MTTKRPs of $\bclY$ with the factor matrices then incur $\clO(Rn^N)$ cost each.
On the other hand, if $\bclY$ is sparse with $\operatorname{nnz}(\bclY)$ nonzeros, then the cost of computing and storing $\bclY$ is reduced to $\clO(\operatorname{nnz}(\bclY))$, and the cost of computing each MTTKRP is reduced to $\clO(r \operatorname{nnz}(\bclY))$ \cite{bader2008efficient}. These are substantial reductions in cost when $\operatorname{nnz}(\bclY) \ll n^N$.

Hence, the idea is to obtain cheap stochastic approximations of the full gradient by randomly subsampling $\bclY$ to produce unbiased sparse approximations $\bctlY$.
We consider two sampling strategies (described in \cref{sec:sampling})
and derive efficient formulas for computing the corresponding stochastic gradients (described in \cref{sec:stochgrad}).

\subsection{Sampling strategies}
\label{sec:sampling}

Here, we briefly review the two sampling strategies (uniform and stratified) from \cite{kolda2020sgf} that we consider.

In uniform sampling, we sample a batch of size $s$, where each sample is an index $i = (i_1, i_2, \dots, i_N)$ drawn uniformly at random from the set of all possible indices $[I_1] \times \cdots \times [I_N]$, with replacement.
Since we are sampling with replacement, an index can be sampled multiple times in a batch;
let $\tls_i$ be the number of times that index $i$ is sampled in a batch.
Now, to make $\bctlY$ an unbiased estimator of $\bclY$, we set
\begin{equation}
    \tly_i = \tls_i \frac{n^N}{s} \frac{\partial \ell}{\partial m_i} (x_i, m_i)
    ,
\end{equation}
for each index $i$ that is sampled.
This sampling strategy is generally most appropriate only if $\bclX$ is dense.

In stratified sampling, we instead sample $p$ nonzero entries and $q$ zero entries,
i.e., we control the number of nonzeros and zeros selected.
Here, to make $\bctlY$ an unbiased estimator of $\bclY$, we set
\begin{equation}
    \tly_i = \tilde{p}_i \frac{\operatorname{nnz}(\bclX)}{p} \frac{\partial \ell}{\partial m_i} (x_i, m_i)
    ,
\end{equation}
for each nonzero entry that is sampled (where $\tlp_i$ is the number of times it is sampled),
and we set
\begin{equation}
    \tly_i = \tilde{q}_i \frac{1 - \operatorname{nnz}(\bclX)}{q} \frac{\partial \ell}{\partial m_i} (x_i, m_i)
    ,
\end{equation}
for each zero entry that is sampled (where $\tlq_i$ is the number of times it is sampled).
This sampling strategy is generally most appropriate when $\bclX$ is sparse; in this case, uniform sampling would rarely sample nonzeros, which can be the most informative entries.
When $\bclX$ is stored in a sparse COOrdinate format, sampling nonzeros can be done efficiently by sampling from the list of stored entries. Sampling zeros can be done in this case via rejection sampling.

\subsection{Computing stochastic gradients}
\label{sec:stochgrad}

With a sparse stochastic approximation $\bctlY$ of $\bclY$ in hand, computing the stochastic gradients can be accomplished by simply replacing $\bclY$ with $\bctlY$ in \cref{eq:symgcp:grad:lam,eq:symgcp:grad:Ak}.
The following theorem shows how to exploit the sparsity of $\bctlY$ to efficiently compute these stochastic gradients.
The approach is essentially the same as for nonsymmetric GCP;
we provide a proof in \cref{sec:stochastic:gradients:proof} for the reader's convenience.

\begin{theorem}[Stochastic SymGCP gradients] \label[theorem]{theorem:stochastic:gradients}
    For a sparse stochastic approximation $\bctlY$ of $\bclY$
    with $s$ nonzero elements at indices $\left(i_1^1, \dots, i_N^1\right), \dots, \left(i_1^s, \dots, i_N^s\right)$,
    the stochastic gradient approximations can be computed as follows
    \begin{align}
        \frac{\partial \clF}{\partial \bmlambda}
        &\approx
        (\ast_{j=1}^N \bhtA_{\sigma_j})^T \bhty
        , \\
        \frac{\partial \clF}{\partial \bmA_k}
        &\approx
        \bhtY_{(n)} (\ast_{j \neq n} \bhtA_{\sigma_j}) \operatorname{diag}(\bmlambda)
        ,
    \end{align}
    where
    \begin{equation*}
        \bhtY_{(n)}
        =
        \begin{bmatrix}
            | & & | \\
            \btlY_{i_1^1,\dots, i_{n-1}^1, :, i_{n+1}^1, \dots, i_N^1} & \!\!\cdots\!\! & \btlY_{i_1^s,\dots, i_{n-1}^s, :, i_{n+1}^s, \dots, i_N^s} \\
            | & & |
        \end{bmatrix}
    \end{equation*}
    and
    \begin{align*}
        \bhtA_{\sigma_j}
        &
        =
        \begin{bmatrix}
            - (\bmA_{\sigma_j})_{i_j^1,:} - \\
            \vdots \\
            - (\bmA_{\sigma_j})_{i_j^s,:} -
        \end{bmatrix}
        , &
        \bhty
        &
        =
        \begin{bmatrix}
            \tly_{i_1^1, \dots, i_N^1} \\
            \vdots \\
            \tly_{i_1^s, \dots, i_N^s}
        \end{bmatrix}
        .
    \end{align*}
\end{theorem}

\section{Numerical Experiments} \label{sec:numerical-experiments}

This section evaluates the proposed SymGCP decomposition through numerical experiments on synthetic datasets.
For nonstochastic fitting of the model, we use the gradients derived in \cref{sec:gradients} with limited-memory BFGS with bound constraints (L-BFGS-B) \cite{Byrd_1995}.
For the stochastic algorithm, we use Adam \cite{kingma2014adam} with standard modifications similar to those used in \cite{kolda2020sgf}. In particular, we separate iterates into epochs with a predefined number of iterations, and estimate the value of the objective function at the end of each epoch. If the value of the objective function fails to decrease by a predefined relative factor $\kappa$ (e.g., $0.99$), we declare it to be a ``bad epoch'' and reset all of the parameters (current factor matrix values, first- and second-order moment estimates, and the iteration counter) to saved values from the previous epoch.

\subsection{Fully symmetric binary tensors} \label{sec:binary-nonstochastic-exp}
Here, we test the effectiveness of SymGCP for recovering the underlying factors of a fully symmetric binary tensor, where the factors use an odds link.
We use the same data generating process as in \cite[Appendix~D]{kolda2020sgf}, modified to make a fully symmetric data tensor. Namely, to make a rank-$r$ $m$-way binary tensor of size $n$ in each mode, we first choose a sparsity level $\delta$, a probability of ones for the signal components $\rho_{\text{high}}$, and a probability of ones for the noise component $\rho_{\text{low}}$. The first $r-1$ columns of the true factor matrix $\bmA_{\ast}$ are designated as signal components, and the last column is designated as a noise component. For each signal component, we randomly select a fraction $\delta$ of the entries to be nonzero, and generate a value from a $\clN (\sqrt[m]{\rho_{\text{high}}/(1 - \rho_{\text{high}})}, 0.5)$ distribution for each nonzero. Every entry in the noise column is set to exactly $\sqrt[m]{\rho_{\text{low}}/(1 - \rho_{\text{low}})}$. Finally, using the true model tensor $\bclM_{\ast}$ created from $\bmA_{\ast}$, each entry $x_i$ in the data tensor $\bclX$ is set to one with probability $m^{\ast}_i/(1 + m^{\ast}_i)$.

To obtain initializations for SymGCP, we generate the initial factor matrix $\bmA$ by first sampling a random Gaussian matrix with $\clN(0,1)$ entries then rescaling it so that the norm of the initialized model tensor $\bclM$ is the same as the norm of the data tensor $\bclX$.
We run SymGCP for each initialization with both the least-squares and Bernoulli (using an odds link) loss functions. To compare how close a recovered solution $\bhtA$ is to the true factor matrix $\bmA_{\ast}$, we compute the cosine similarity score
\begin{equation}
    \frac{1}{r}\sum_{j=1}^r \cos((\bmA_{\ast})_{:,j}, (\bhtA)_{:,\pi(j)}),
\end{equation}
where $\pi$ is the permutation of the columns of $\bhtA$ that produces the highest score.

For our experiments, we set $m=4$, $n=50$, $r=5$, $\rho_{\text{high}}=0.9$, $\rho_{\text{low}} = 0.002$, and $\delta=0.15$. We created one true factor matrix $\bmA_{\ast}$ then generated $20$ different instances of the data tensor $\bclX$ from $\bmA_{\ast}$. For each instance, we ran $25$ different initializations of SymGCP for both losses.

\begin{figure}
    \centering
    \includegraphics[width=0.95\linewidth]{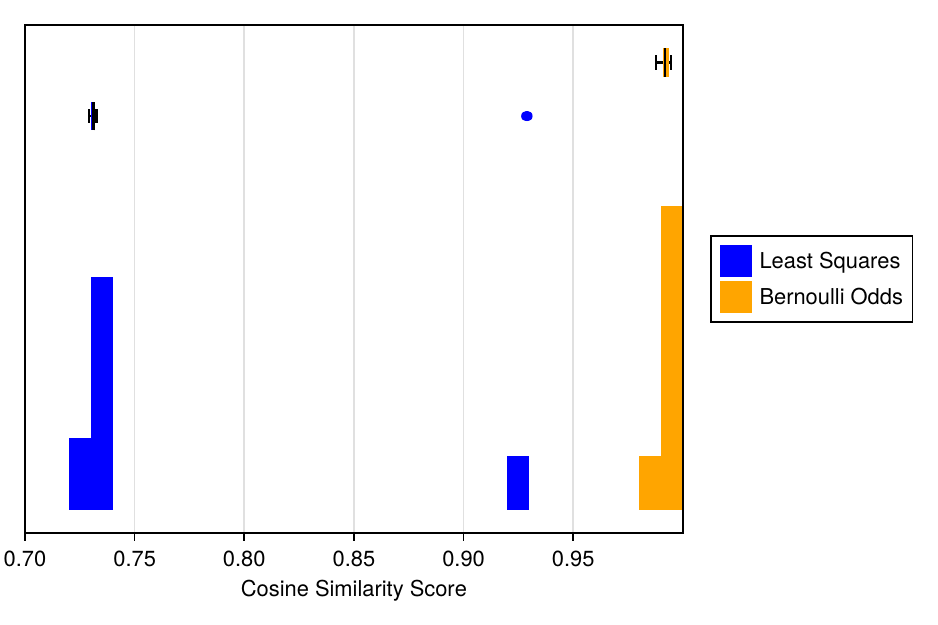}
    \caption{Histograms of best cosine similarity scores for each instance for SymGCP on binary test problems with the least-squares and Bernoulli odds loss functions. The best similarity score for each instance is the similarity score for the initialization which achieved the lowest final loss value. The corresponding boxplots are shown above the histograms.}
    \label{fig:binary-results}
\end{figure}

\Cref{fig:binary-results} compares the best cosine similarity scores for each instance for the least-squares and Bernoulli odds losses. SymGCP with the Bernoulli odds loss better recovers the underlying factors in $\bmA_{\ast}$, with a median best score of $0.998$ across the $20$ different instances, compared to a median best score of $0.835$ for least-squares. Furthermore, the worst score for Bernoulli loss across instances ($0.997$) is better than the best score for least-squares loss across instances ($0.982$).
This illustrates the benefit of using general losses to model the noise statistics.

\subsection{Stochastic experiments}
We next test the stochastic approach using the stochastic gradients from \cref{sec:stochastic-gradients} for SymGCP on simulated fully symmetric binary tensor data generated in the same way as in \cref{sec:binary-nonstochastic-exp}.
Using the same values for all of the experimental parameters, we created a new true factor matrix $\bmA_{\ast}$ and then again generated $20$ different instances of the data tensor $\bclX$ from $\bmA_{\ast}$. The data tensors are between $0.43\%$ and $0.44\%$ percent sparse, with a total of between $26874$ and $27272$ nonzeros.

We run SymGCP with Bernoulli (odds link) loss using L-BFGS-B and Adam with both uniform and stratified sampling for $25$ different initializations for each instance. For Adam, we use a total sample size of $1000$ for both sampling strategies, with $500$ zeros and $500$ non-zeros for stratified.

\begin{figure}
    \centering
    \includegraphics[width=0.98\linewidth]{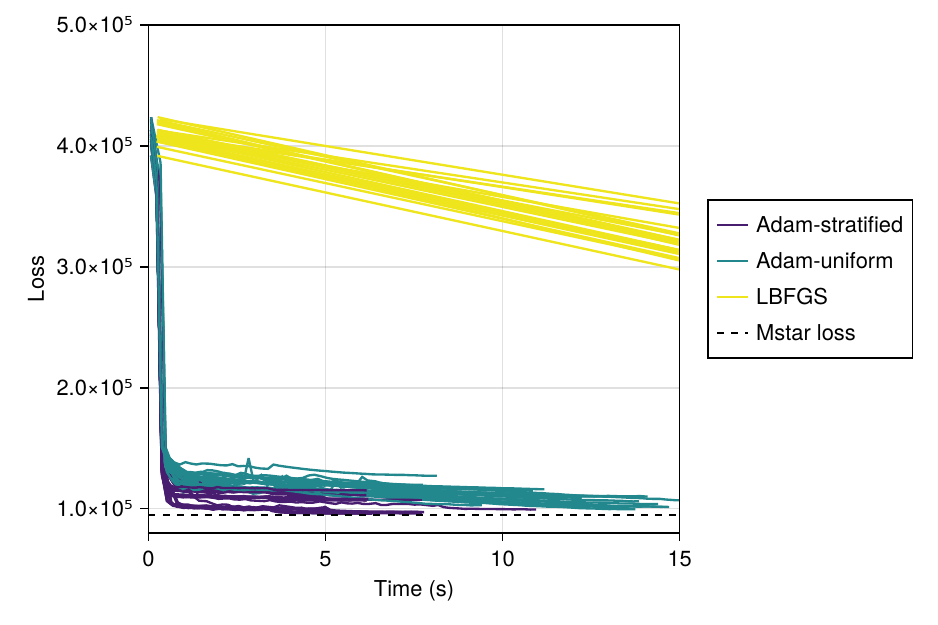}
    \includegraphics[width=0.98\linewidth]{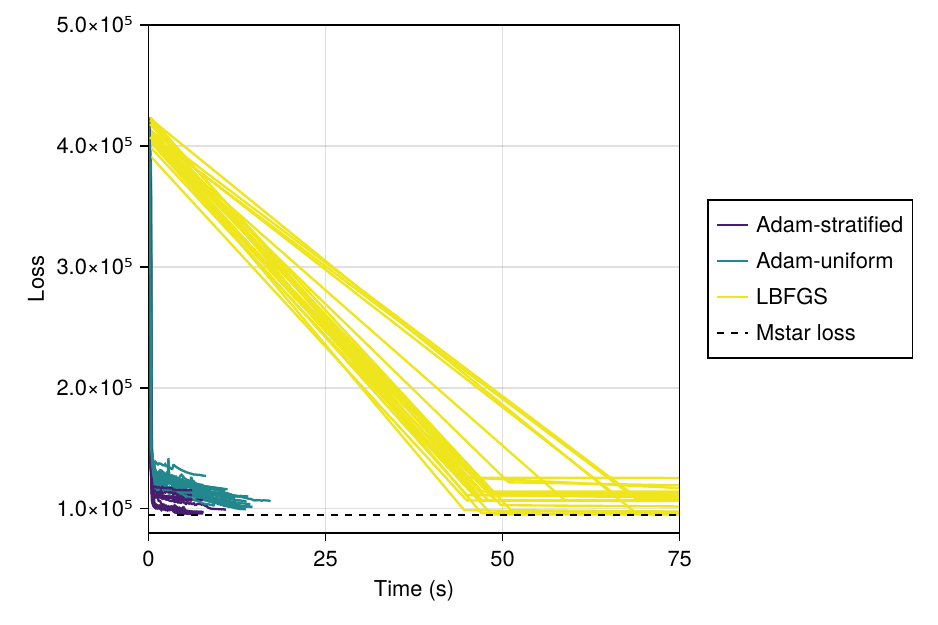}
    \includegraphics[width=0.95\linewidth]{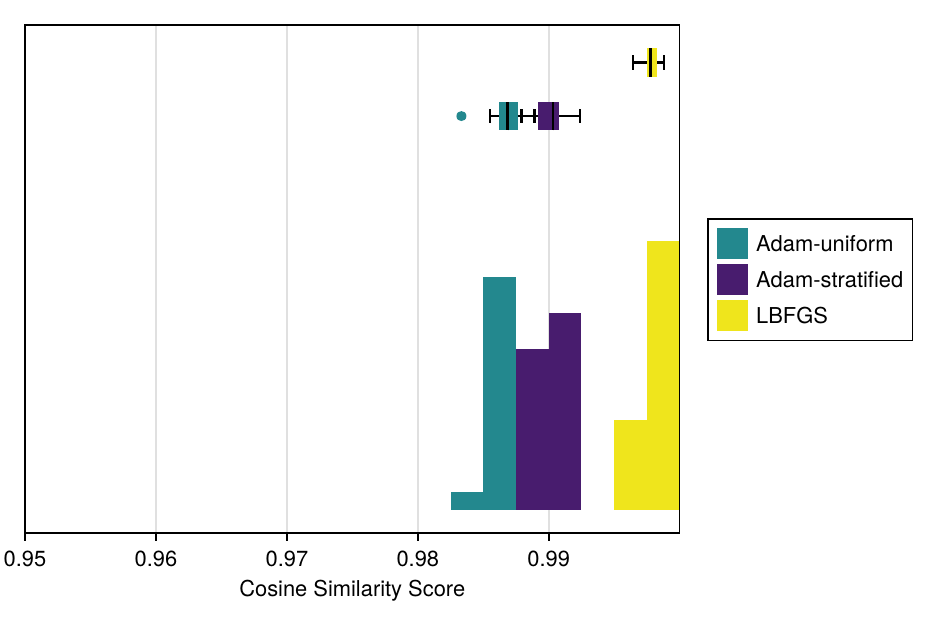}
    \caption{Top: Objective function value over time for SymGCP on a binary test problem for a representative data instance with Adam (stratified and uniform sampling) and L-BFGS-B. Middle: Objective function value over time with zoomed out time-axis. Bottom: Histogram with boxplot showing the best cosine similarity scores for each instance for SymGCP on binary test problems with Adam and L-BFGS-B. The best similarity score for each instance is the similarity score for the initialization that achieved the lowest final loss value.}
    \label{fig:stochastic-results}
\end{figure}

\Cref{fig:stochastic-results} shows the objective function value versus walltime for SymGCP with Adam (uniform and stratified sampling) and L-BFGS-B for a representative data instance (other instances are similar) and the cosine similarity scores for the best initialization for each method across all data instances. Using Adam with stratified sampling for these sparse symmetric binary tensors enables a roughly 50 times speed-up, with most of the L-BFGS-B runs taking around 50 seconds to reach a loss value near the ``true loss'' given by $\bclM_{\ast}$, while most of the stratified Adam runs take around 1 second to reach the same level. Adam with uniform sampling is also much faster than L-BFGS-B, although it is slower than stratified sampling. Overall, the stochastic methods have much cheaper iterations and are thus able to more rapidly descend down the objective function. On the other hand, the stochastic methods also tended to have a lower accuracy for the final solution since they do not converge to as good a solution. Nevertheless, the stochastic solutions were still fairly close to the true underlying factors, with a median cosine similarity score of $0.990$ for stratified sampling and $0.987$ for uniform sampling, as compared to a median score of $0.998$ for L-BFGS-B.

\section{Real Data Experiments} \label{sec:real-data-experiments}

This section evaluates the proposed SymGCP decomposition on two real datasets.
We use the same algorithms as in \cref{sec:numerical-experiments}.

\subsection{Monkey neural data}
We first apply SymGCP to neural data recorded during experiments where a monkey attempted to move a cursor to one of four different targets using a brain machine interface (BMI). This data was original collected to study neural dynamics during a motor learning task \cite{vyas2018neural, vyas2020causal}. We use the data files organized for a demo of the Tensor Toolbox package \cite{kolda2021bmi}. The dataset contains a tensor of size $43 \times 200 \times 88$ which contains the activity for 43 different neurons sampled at 200 consecutive time steps for 88 different trials, where in each trial the monkey attempted to move the cursor to one of the four targets. Target labels for each trial are also provided.

Instead of decomposing the $43 \times 200 \times 88$ neural activity tensor, which we denote by $\bclX$, we consider a tensor comprised of neural coactivations. The $(i,j,k)$ element of the coactivations tensor $\bclC$ is the inner product of the mode-2 fibers of $\bclX$ corresponding to the time-series for neurons $i$ and $j$ during trial $k$:
\begin{equation}
    c_{i,j,k} = (x_{i,:,k})^T x_{j,:,k},
\end{equation}
for all $i \in [43], j \in [43], k \in [88]$. Note that $\bclC$ is symmetric along its first two modes.

We compute SymGCP decompositions of the coactivations tensor $\bclC$ with a nonnegative least-squares loss and varying ranks from $[2,4,6,8,10,15,20]$. We use the target labels for the trials to evaluate how well the decompositions recover the underlying trial structure. To evaluate a rank-$r$ decomposition, we treat each row in the recovered trial factor matrix as a data point in $\bbR^r$, and use K-means with 4 centroids to cluster the rows. After clustering, we record the percentage of trials (rows of the trial factor matrix) that were correctly clustered (after permuting the arbitrary labels of the centroids to have the highest accuracy).

\begin{figure}
    \centering
    \includegraphics[width=0.95\linewidth]{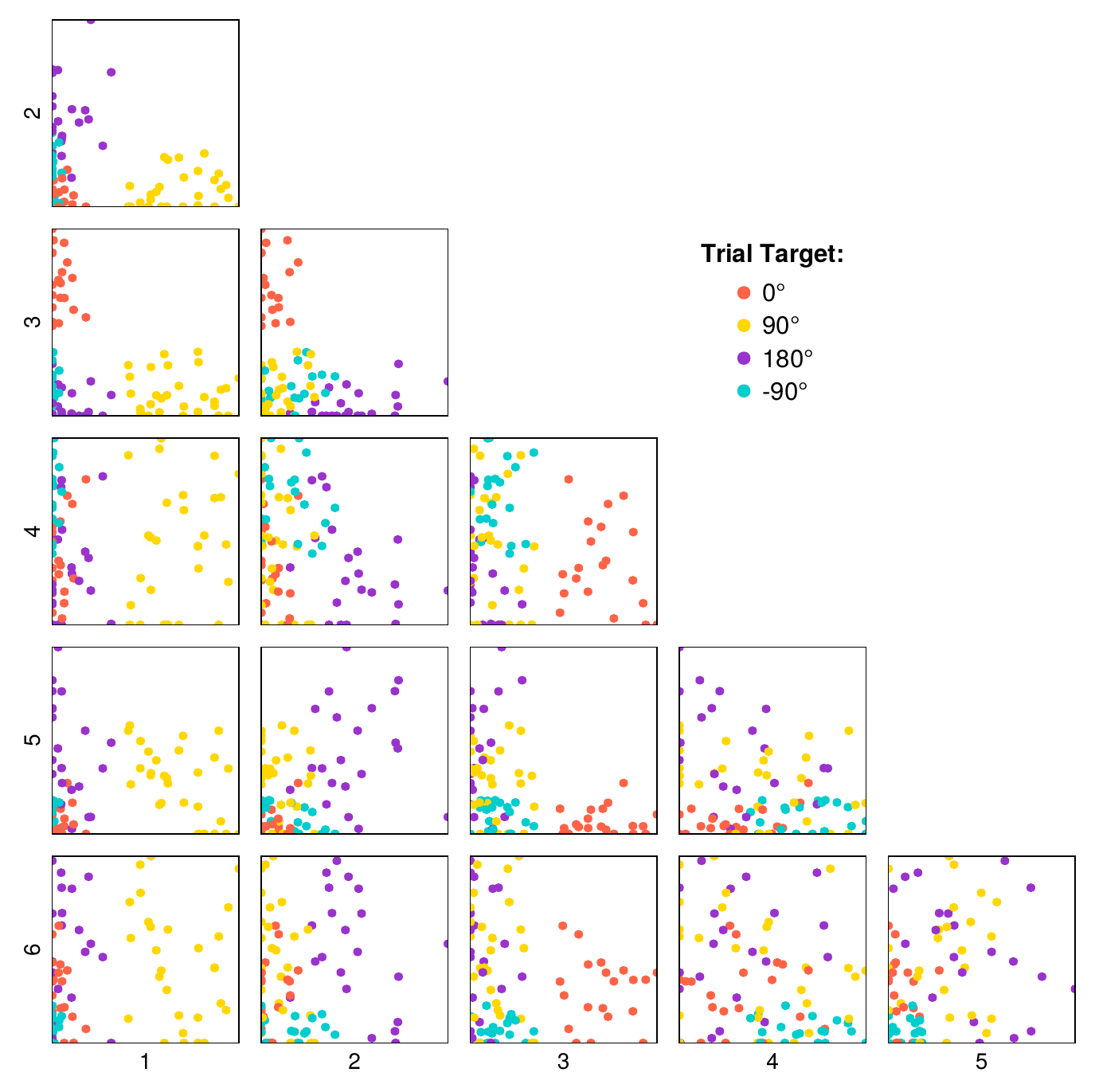}
    \caption{Pairwise relationships between rows of the recovered trial factor matrix for a rank-6 symmetric nonnegative least-squares decomposition of the coactivation tensor from monkey neural data. Trials are colored by the target condition.}
    \label{fig:monkey-pairplot}
\end{figure}

We found that SymGCP with a small rank was able to recover the trial structure of the data, with perfect clustering after applying K-means on the recovered factor matrix for a rank of $6$. In \cref{fig:monkey-pairplot}, we plot the pairwise relationships between the 88 rows (corresponding to the 88 trials) of the recovered factor matrix for a rank-3 decomposition, with the points colored according to which of the 4 target conditions they correspond to. Even though the coactivations tensor contains no explicit information about the target condition for each trial, the trial factor matrix from the symmetric decomposition clearly distinguishes between the different target conditions.

\subsection{UC Irvine social network}
\begin{figure}
    \centering
    \includegraphics[width=0.95\linewidth]{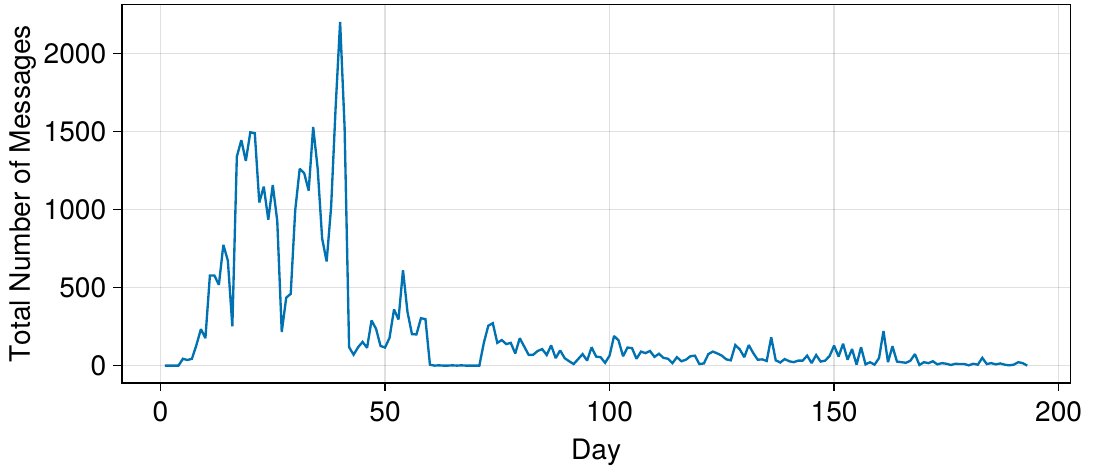}
    \caption{Total number of messages exchanged between the $200$ most frequent users in UC Irvine social network dataset \cite{opsahl2009clustering} on each day of the experiment.}
    \label{fig:ucirvine-time}
\end{figure}

Next, we evaluate SymGCP on the UC Irvine social network dataset \cite{opsahl2009clustering}, which contains a daily count of the number of messages between pairs of users of a social network, with a total of $1899$ users and $193$ days where messages were sent. Starting with the original dataset which contains the number of messages between each pair of users on each day, we construct a count tensor where element $(x,y,z)$ contains the number of messages exchanged between users $x$ and $y$ on day $z$, and is therefore symmetric along the first two modes. We include only the $200$ most active users. While we do not have any information about underlying structure in the data, we can observe that there is significantly more overall communication earlier on (roughly from day $0$ to day $50$) in the experiment, as illustrated in \cref{fig:ucirvine-time}, which shows the total number of messages between all of the top $200$ users for each day.

\begin{figure}
    \centering
    \includegraphics[width=0.49\linewidth]{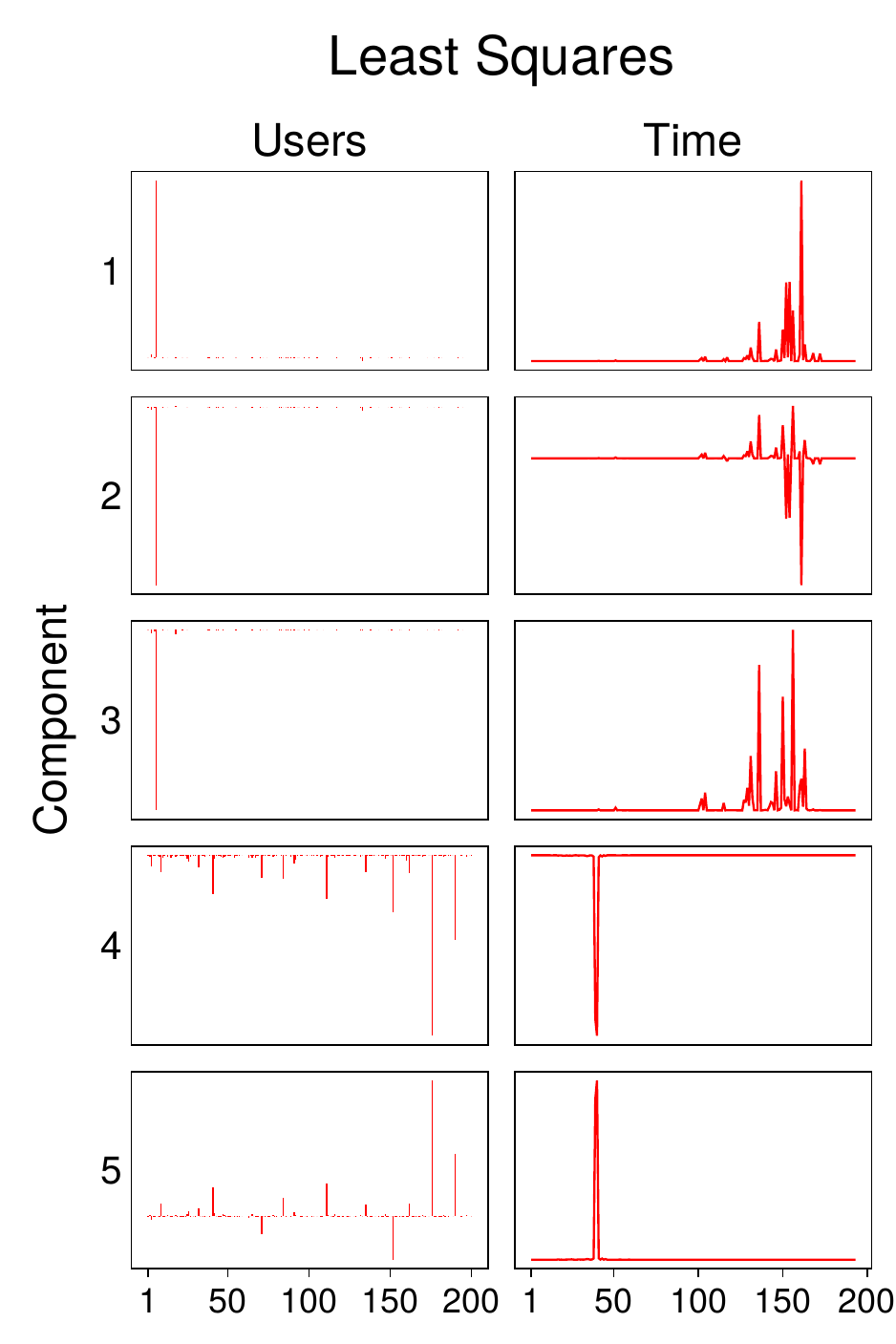}
    \includegraphics[width=0.49\linewidth]{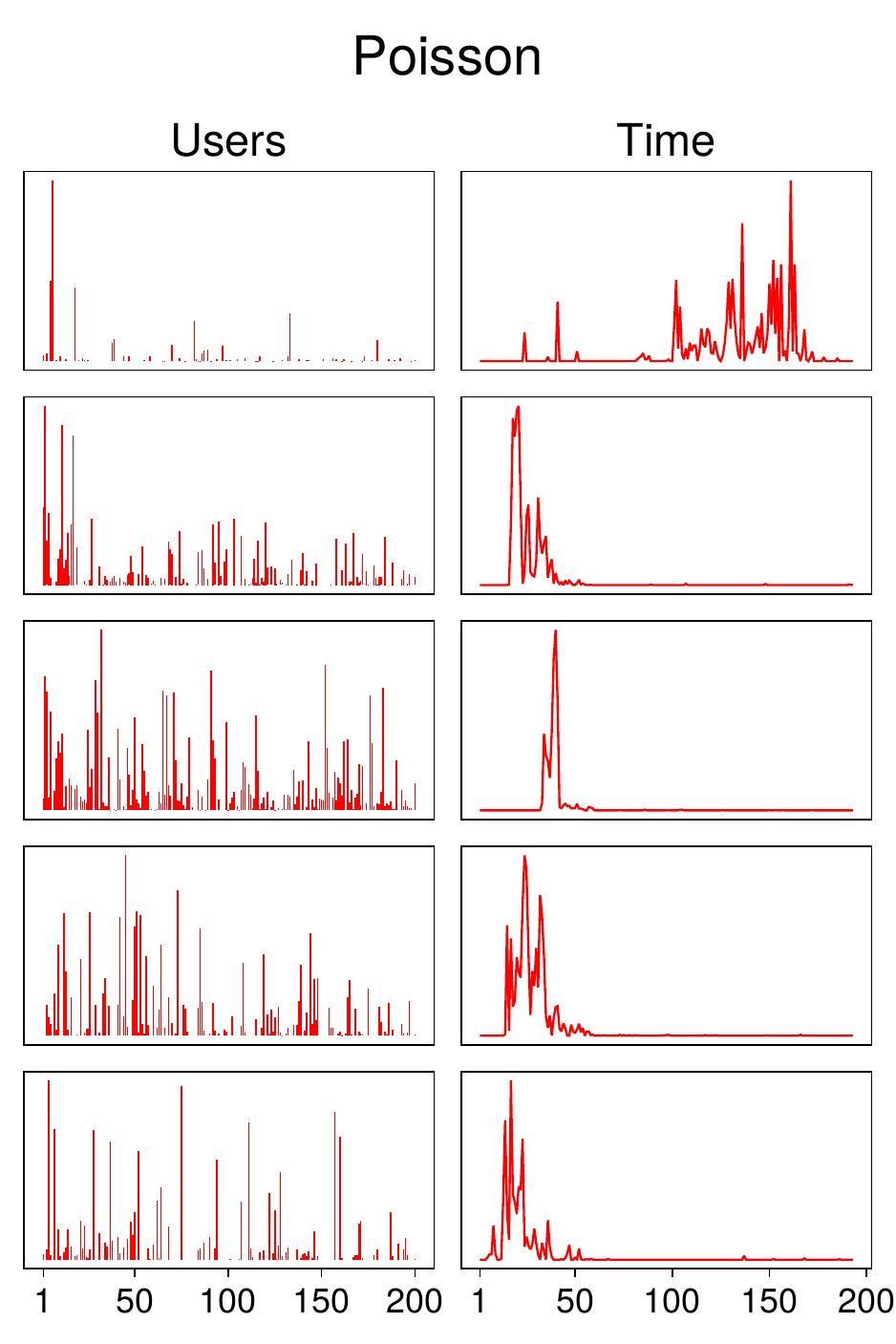}
    \caption{SymGCP decompositions with a least-squares loss (left, corresponding to symmetric CP) and a Poisson loss (right) for the top $200$ most active users in the UC Irvine social network dataset \cite{opsahl2009clustering}. The top $5$ components (in terms of largest weight in $\lambda$) are shown for a rank $r=10$ decomposition for each.}
    \label{fig:ucirvine-results}
\end{figure}

For the constructed count tensor, we compute a rank $r=10$ symmetric CP decomposition (least-squares loss) and a SymGCP decomposition with a loss function corresponding to a Poisson distribution. \Cref{fig:ucirvine-results} shows the results for the top $5$ components (after sorting by descending weight) for both decompositions. As opposed to symmetric CP decomposition, which recovers components that focus only on a few users and days (the top three components appear to all identify the same single user, whereas components 4 and 5 identify a single day), SymGCP with the Poisson loss discovers components which include more users and longer time intervals, likely indicating more meaningful social groups. Components 1,2,4, and 5 specifically identify groups of users with high activity early in the experiment, matching the overall high activity early on in \cref{fig:ucirvine-time}.

\section{Conclusion} \label{sec:conclusion}
This paper developed a new CP decomposition that can be used for tensors with arbitrary forms of symmetry and with general loss functions. We defined a general notion of tensor symmetry that includes the fully symmetric tensors that have been the primary focus of previous symmetric CP works \cite{kolda2015nof, sherman2020estimating, kileel2019subspace}, as well as partial symmetry. We derived the gradients of the SymGCP optimization problem and applied an all-at-once gradient-based method. We also developed efficient stochastic approximations of the gradients that enable scalability to very large tensors. Numerical experiments demonstrated that SymGCP with appropriately chosen loss functions can effectively recover the underlying low-rank factors in symmetric binary tensors, with speed-ups available for large tensors via the stochastic approach. Finally, we applied SymGCP to monkey neural data \cite{vyas2018neural,vyas2020causal,kolda2021bmi} and social network data \cite{opsahl2009clustering} to find meaningful underlying components in the data.

One possible future direction of work is developing more efficient optimization algorithms for computing SymGCP decompositions, potentially by further exploiting symmetric structure within the data tensor and/or model. For example, the intermediate results from the sequences of MTTKRPs for one factor matrix may be able to be used to reduce the computation needed for the MTTKRPs of the other factor matrices. Another future direction is developing streaming/online algorithms for SymGCP that could be used to compute symmetric generalized CP decompositions of data as it is generated, e.g., in settings such as large-scale simulations.

\appendices
\crefalias{section}{appendix}

\section{Proof of Proposition~\ref{thm:symmetric:ktensor:issymmetric}} \label{sec:symmetric:ktensor:issymmetric}

Let $i \in [I_1] \times \cdots \times [I_N]$ and $\pi \in \Pi(\clI_1) \times \cdots \times \Pi(\clI_K)$ be arbitrary. Note that
\begin{multline}
    \llbracket\bmlambda ; \clI_1 \to \bmA_{1}, \dots, \clI_K \to \bmA_{K}\rrbracket_i
    \\
    \begin{aligned}
        &=
        \llbracket\bmlambda ; \bmA_{\sigma_1}, \dots, \bmA_{\sigma_N} \rrbracket_i
        \\
        &=
        \left( \sum_{j=1}^r \lambda_j \cdot \left(\bmA_{\sigma_1}\right)_{:,j} \circ \cdots \circ \left(\bmA_{\sigma_N}\right)_{:,j} \right)_i
        \\
        &=
        \sum_{j=1}^r \lambda_j \cdot \left(\bmA_{\sigma_1}\right)_{i_1,j} \circ \cdots \circ \left(\bmA_{\sigma_N}\right)_{i_N,j}
        \\
        &=
        \sum_{j=1}^r \lambda_j \prod_{n=1}^N \left(\bmA_{\sigma_N}\right)_{i_n,j}
        \\
        &=
        \sum_{j=1}^r \lambda_j \prod_{\ell=1}^K \prod_{k \in \clI_{\ell}} \left(\bmA_{\ell}\right)_{i_k,j}.
    \end{aligned}
\end{multline}
Now, since $\pi$ only permutes indices within the cells of the partition, $\{ \pi(k): k \in \clI_{\ell} \} = \clI_{\ell}$ for any $\ell \in [K]$, and so
\begin{equation}
    \prod_{k \in \clI_{\ell}} \left(\bmA_{\ell}\right)_{i_{\pi(k)},j}
    =
    \prod_{k \in \clI_{\ell}} \left(\bmA_{\ell}\right)_{i_k,j}
\end{equation}
for any $j \in [r]$ and $\ell \in [K]$. Therefore
\begin{multline}
    \llbracket\bmlambda ; \clI_1 \to \bmA_{1}, \dots, \clI_K \to \bmA_{K}\rrbracket_{\pi(i)}
    \\
    \begin{aligned}
        &=
        \sum_{j=1}^r \lambda_j \prod_{\ell=1}^K \prod_{k \in \clI_{\ell}} \left(\bmA_{\ell}\right)_{i_{\pi(k)},j}
        \\
        &=
        \sum_{j=1}^r \lambda_j \prod_{\ell=1}^K \prod_{k \in \clI_{\ell}} \left(\bmA_{\ell}\right)_{i_k,j}
        \\
        &=
        \llbracket\bmlambda ; \clI_1 \to \bmA_{1}, \dots, \clI_K \to \bmA_{K}\rrbracket_{i},
    \end{aligned}
\end{multline}
which concludes the proof.
\qed

\section{Proof of \cref{theorem:symgcp-gradients}}
\label{theorem:symgcp-gradients:proof}

Recall from \cite[Corollary 4]{hong2020gcp} that the nonsymmetric GCP objective function
\begin{align}
    \clG(\bmlambda, \btlA_{1}, \dots, \btlA_{N})
    &= \clL (\bclX, \llbracket\bmlambda ; \btlA_{1}, \dots, \btlA_{N}\rrbracket)
\end{align}
has the gradients
\begin{align}
    \frac{\partial \clG}{\partial \bmlambda}
    &=
    (\odot_{j=1}^N \btlA_{j})^T \operatorname{vec}(\bclY),
    \\
    \frac{\partial \clG}{\partial \btlA_k}
    &=
    \bmY_{(k)} (\odot_{j \neq k} \btlA_{j}) \operatorname{diag} (\bmlambda),
\end{align}
where
\begin{align*}
    \odot_{j=1}^N \btlA_{j}
    &=
    \btlA_{N} \odot \cdots \odot \btlA_{1}
    , \\
    \odot_{j \neq k} \btlA_{j}
    &=
    \btlA_{N} \odot \cdots \odot \btlA_{k+1} \odot \btlA_{k-1} \odot \cdots \odot \btlA_{1}
    ,
\end{align*}
and $\bclY$ is the derivative tensor whose entries are
\begin{equation*}
    y_i = \frac{\partial \ell}{\partial m_i} (x_i, m_i)
    .
\end{equation*}
For SymGCP, we now have a factor matrix for each cell of symmetric modes rather than each mode.
Note next that
\begin{align}
    &
    \clF\left(\bmlambda, \bmA_{1}, \dots, \bmA_{K}\right)
    \\&\qquad
    =
    \clL \left(\bclX, \llbracket\bmlambda ; \clI_1 \to \bmA_{1}, \dots, \clI_K \to \bmA_{K}\rrbracket\right)
    \nonumber
    \\&\qquad
    =
    \clL \left(\bclX, \llbracket\bmlambda ; \bmA_{\sigma_1}, \dots, \bmA_{\sigma_N}\rrbracket\right)
    \nonumber
    \\&\qquad
    =
    \clG\left(\bmlambda, \bmA_{\sigma_1}, \dots, \bmA_{\sigma_N}\right)
    ,
\end{align}
and so it immediately follows that
\begin{equation}
    \frac{\partial \clF}{\partial \bmlambda}
    =
    \left(\odot_{j=1}^N \bmA_{\sigma_j}\right)^T \operatorname{vec}\left(\mathcal{Y}\right).
\end{equation}
For gradients with respect to the factor matrices,
applying chain rule yields
\begin{align}
    \frac{\partial \clF}{\partial \operatorname{vec}(\bmA_k)}
    &=
    \sum_{t=1}^N \left[ \frac{\partial \operatorname{vec} (\bmA_{\sigma_t})}{\partial \operatorname{vec} (\bmA_k)} \right]^T \frac{\partial \clG}{\partial \operatorname{vec}(\bmA_{\sigma_t})}
    .
\end{align}
Now, for each $t \in [N]$, we have
\begin{equation*}
    \frac{\partial \operatorname{vec} (\bmA_{\sigma_t})}{\partial \operatorname{vec} (\bmA_k)}
    =
    \begin{cases}
        \bmI , & \text{if } \sigma_t=k
        , \\
        \bm0 , & \text{otherwise}
        ,
    \end{cases}
\end{equation*}
so we finally have
\begin{equation}
    \frac{\partial \clF}{\partial \bmA_k}
    =
    \sum_{t \in \clI_k} \frac{\partial \clG}{\partial \bmA_{\sigma_t}}
    =
    \sum_{t \in \clI_k} \bmY_{(t)} \left(\odot_{j \neq t} \bmA_{\sigma_j}\right) \operatorname{diag} \left(\bmlambda\right)
    ,
\end{equation}
which completes the proof.
\qed

\section{Proof of Proposition~\ref{lem:symmetric-mttkrps}}
\label{lem:symmetric-mttkrps:proof}

Assume without loss of generality that $u < v$. Using well-known properties of Kruskal tensors, we can write the inner product of $\bmY_{(u)}\left(\odot_{j \neq u} \bmA_{\sigma_j}\right)$ with an arbitrary matrix $\bmB \in \bbR^{s \times r}$, where $s = I_u = I_v$, as
\begin{multline}
    \left\langle \bmB,  \bmY_{(u)}\left(\odot_{j \neq u} \bmA_{\sigma_j}\right) \right\rangle
    \\
    \begin{aligned}
        &=
        \left\langle \bmY_{(u)}, \bmB \left(\odot_{j \neq u} \bmA_{\sigma_j}\right)^T \right\rangle
        \\
        &=
        \left\langle \mathcal{Y}, \llbracket \bmA_{\sigma_1}, \dots, \bmA_{\sigma_{u-1}}, \bmB, \bmA_{\sigma_{u+1}}, \dots, \bmA_{\sigma_N} \rrbracket \right\rangle.
    \end{aligned}
\end{multline}
Now, since $\bclY$ is symmetric along modes $u$ and $v$, and $\bmA_{\sigma_u} = \bmA_{\sigma_v}$, we have that
\begin{align}
    \Big\langle \mathcal{Y}, &\llbracket \bmA_{\sigma_1}, \dots, \bmA_{\sigma_{u-1}}, \bmB, \bmA_{\sigma_{u+1}}, \dots, \bmA_{\sigma_N} \rrbracket \Big\rangle
    \notag
    \\
    &=
    \Big\langle \mathcal{Y}, \llbracket \bmA_{\sigma_1}, \dots, \bmA_{\sigma_{u-1}}, \bmA_{\sigma_v}, \bmA_{\sigma_{u+1}},
    \\&\nonumber \qquad \qquad
    \dots,
    \bmA_{\sigma_{v-1}}, \bmB, \bmA_{\sigma_{v+1}}, \dots, \bmA_{\sigma_N} \rrbracket \Big\rangle
    \notag
    \\
    &=
    \Big\langle \mathcal{Y}, \llbracket \bmA_{\sigma_1}, \dots, \bmA_{\sigma_{u-1}}, \bmA_{\sigma_u}, \bmA_{\sigma_{u+1}},
    \\&\nonumber \qquad \qquad
    \dots,
    \bmA_{\sigma_{v-1}}, \bmB, \bmA_{\sigma_{v+1}}, \dots, \bmA_{\sigma_N} \rrbracket \Big\rangle
    \notag
    \\
    &=
    \left\langle \bmY_{(v)}, \bmB \left(\odot_{j \neq v} \bmA_{\sigma_j}\right)^T \right\rangle
    \notag
    \\
    &=
    \left\langle \bmB,  \bmY_{(v)}\left(\odot_{j \neq v} \bmA_{\sigma_j}\right) \right\rangle
    .
\end{align}
Thus, we have shown that for any $\bmB \in \bbR^{s \times r}$
\begin{equation*}
    \left\langle \bmB,  \bmY_{(u)}\left(\odot_{j \neq u} \bmA_{\sigma_j}\right) \right\rangle
    =
    \left\langle \bmB,  \bmY_{(v)}\left(\odot_{j \neq v} \bmA_{\sigma_j}\right) \right\rangle
    ,
\end{equation*}
which implies that
$\bmY_{(u)}\left(\odot_{j \neq u} \bmA_{\sigma_j}\right)
=
\bmY_{(v)}\left(\odot_{j \neq v} \bmA_{\sigma_j}\right)$,
thus completing the proof.
\qed

\section{Proof of \cref{theorem:stochastic:gradients}} \label{sec:stochastic:gradients:proof}

We first show that $\btlY_{(n)} \left(\odot_{j \neq n} \bmA_{\sigma_j}\right) = \bhtY_{(n)} \left(\ast_{j \neq n} \bhtA_{\sigma_j}\right)$. We start by expanding out the Khatri-Rao product into a columnwise Kronecker product and then applying the definition of matrix multiplication. If we let $\left(\otimes_{j\neq n} \left(\bmA_{\sigma_j}\right)_{:,\ell}\right)$ be a shorthand for
\begin{equation*}
    \left(\bmA_{\sigma_N}\right)_{:,\ell} \otimes \cdots \otimes \left(\bmA_{\sigma_{n+1}}\right)_{:,\ell} \otimes \left(\bmA_{\sigma_{n-1}}\right)_{:,\ell} \otimes \cdots \otimes \left(\bmA_{\sigma_1}\right)_{:,\ell},
\end{equation*}
and $\sum_{\substack{i_k = 1 \\ k \neq n}}^{I_k}$ be a shorthand for
\begin{equation*}
    \sum_{i_1 = 1}^{I_1} \cdots \sum_{i_{n-1} = 1}^{I_{n-1}} \; \sum_{i_{n+1} = 1}^{I_{n+1}} \cdots \sum_{i_N = 1}^{I_N},
\end{equation*}
then the $(u,v)$ element of the mode-$n$ MTTKRP $\btlY_{(n)} \left(\odot_{j \neq n} \bmA_{\sigma_j}\right)$ is
\begin{multline}
    \left[ \btlY_{(n)} \left(\odot_{j \neq n} \bmA_{\sigma_j}\right) \right]_{u,v}
    \\
    \begin{aligned}
        &=
        \left[ \btlY_{(n)}
            \left[
                \begin{array}{@{}ccc@{}}
                    \left(\otimes_{j\neq n} \left(\bmA_{\sigma_j}\right)_{:,1}\right) & \cdots & \left(\otimes_{j\neq n} \left(\bmA_{\sigma_j}\right)_{:,r}\right)
                \end{array}
            \right]
        \right]_{u,v}
        \\
        &=
        \sum_{\substack{i_k = 1 \\ k \neq n}}^{I_k}
        \left( \tilde{y}_{i_1,\dots,i_{n-1},u,i_{n+1},\dots,i_N} \left(\prod_{j \neq n} \left(\bmA_{\sigma_j}\right)_{i_j,v}\right) \right)
    \end{aligned}
\end{multline}
Since $\bctlY$ is sparse with only $s$ nonzeros, the above sum reduces to
\begin{multline}
    \left[\btlY_{(n)} \left(\odot_{j \neq n} \bmA_{\sigma_j}\right)\right]_{u,v}
    \\
    =
    \sum_{w=1}^s \left( \tilde{y}_{i_1^w,\dots,i_{n-1}^w,u,i_{n+1}^w,\dots,i_N^w} \left(\prod_{j \neq n} \left(\bmA_{\sigma_j}\right)_{i^w_j,v}\right) \right)
\end{multline}
Now, if we define $\bhtY_{(n)}$ as
\begin{multline}
    \bhtY_{(n)}
    =
    \\
    \left[
        \begin{array}{@{}ccc@{}}
            \tly_{i_1^1,\dots, i_{n-1}^1, 1, i_{n+1}^1, \dots, i_N^1} & \dots & \tly_{i_1^s,\dots, i_{n-1}^s, 1, i_{n+1}^s, \dots, i_N^s}
            \\
            \vdots & \vdots & \vdots
            \\
            \tly_{i_1^1,\dots, i_{n-1}^1, I_n, i_{n+1}^1, \dots, i_N^1} & \dots & \tly_{i_1^s,\dots, i_{n-1}^s, I_n, i_{n+1}^s, \dots, i_N^s}
        \end{array}
    \right],
\end{multline}
and $\bhtA_{\sigma_j}$ as
\begin{equation}
    \bhtA_{\sigma_j}
    =
    \left[
        \begin{array}{@{}c@{}}
            - (\bmA_{\sigma_j})_{i_j^1,:} - \\
            \vdots \\
            - (\bmA_{\sigma_j})_{i_j^s,:} -
        \end{array}
    \right]
    =
    \left[
        \begin{array}{@{}ccc@{}}
            (\bmA_{\sigma_j})_{i_j^1, 1} & \cdots & (\bmA_{\sigma_j})_{i_j^1, r} \\
            \vdots & \vdots & \vdots \\
            (\bmA_{\sigma_j})_{i_j^s, 1} & \cdots & (\bmA_{\sigma_j})_{i_j^s, r}
        \end{array}
    \right],
\end{equation}
then we have that
\begin{multline}
    \left[ \bhtY_{(n)} \left(\ast_{j \neq n} \bhtA_{\sigma_j}\right)\right]_{u,v}
    \\
    \begin{aligned}
        &=
        \left[\bhtY_{(n)}
            \left[
                \begin{array}{@{}ccc@{}}
                    \prod_{j \neq n} \left(\bmA_{\sigma_j}\right)_{i^1_j,1} & \cdots & \prod_{j \neq n} \left(\bmA_{\sigma_j}\right)_{i^1_j,r} \\
                    \vdots & \vdots & \vdots \\
                    \prod_{j \neq n} \left(\bmA_{\sigma_j}\right)_{i^s_j,1} & \cdots & \prod_{j \neq n} \left(\bmA_{\sigma_j}\right)_{i^s_j,r}
                \end{array}
            \right]
        \right]_{u,v}
        \\
        &=
        \sum_{w=1}^s \left( \tilde{y}_{i_1^w,\dots,i_{n-1}^w,u,i_{n+1}^w,\dots,i_N^w} \left(\prod_{j \neq n} \left(\bmA_{\sigma_j}\right)_{i^w_j,v}\right) \right)
        \\
        &=
        \left[ \btlY_{(n)} \left(\odot_{j \neq n} \bmA_{\sigma_j}\right)\right]_{u,v}.
    \end{aligned}
\end{multline}
Next, we will show that $\left(\odot_{j=1}^N \bmA_{\sigma_j}\right)^T \operatorname{vec}(\bctlY) = \bhty\left(\ast_{j=1}^N \bhtA_{\sigma_j}\right)$ using a similar approach. If we let
\begin{equation}
    \left(\otimes_{j=1}^n \left(\bmA_{\sigma_j}\right)_{:,\ell}\right)
    =
    \left(\bmA_N\right)_{:,\ell} \otimes \cdots \otimes \left(\bmA_1\right)_{:,\ell},
\end{equation}
then the $v$th element of $(\odot_{j=1}^N \bmA_{\sigma_j})^T \operatorname{vec}(\bctlY)$ is
\begin{multline}
    \left[ \left(\odot_{j=1}^N \bmA_{\sigma_j}\right)^T \operatorname{vec}\left(\bctlY\right)\right]_v
    \\
    \begin{aligned}
        &=
        \left[
            \left[
                \begin{array}{@{}ccc@{}}
                    \left(\otimes_{j=1}^n \left(\bmA_{\sigma_j}\right)_{:,1}\right) & \cdots & \left(\otimes_{j=1}^n \left(\bmA_{\sigma_j}\right)_{:,r}\right)
                \end{array}
            \right]
        \operatorname{vec}\left(\bctlY\right)\right]_v
        \\
        &= \sum_{i_1 = 1}^{I_1} \cdots \sum_{i_N = 1}^{I_N} \left( \tilde{y}_{i_1,\dots,i_N} \left(\prod_{j=1}^N \left(\bmA_{\sigma_j}\right)_{i_j,v}\right) \right)
    \end{aligned}
\end{multline}
Since $\bctlY$ is sparse with only $s$ nonzero elements, the above sum reduces to
\begin{multline}
    \left[ \left(\odot_{j=1}^N \bmA_{\sigma_j}\right)^T \operatorname{vec}\left(\bctlY\right)\right]_v
    \\
    =
    \sum_{w=1}^s \left( \tilde{y}_{i_1^w,\dots,i_N^w} \left(\prod_{j=1}^N \left(\bmA_{\sigma_j}\right)_{i^w_j,v}\right) \right)
\end{multline}
Now, letting
\begin{equation}
    \bhty =
    \begin{bmatrix}
        \tilde{y}_{i_1^1,\dots,i_N^1} \\
        \vdots \\
        \tilde{y}_{i_1^s,\dots,i_N^s}
    \end{bmatrix}^T,
\end{equation}
and defining $\bhtA_{\sigma_j}$ as before, we have that
\begin{multline}
    \left[ \bhty\left(\ast_{j=1}^N \bhtA_{\sigma_j}\right)\right]_v
    \\
    \begin{aligned}
        &=
        \left[ \bhty
            \left[
                \begin{array}{@{}ccc@{}}
                    \prod_{j=1}^N \left(\bmA_{\sigma_j}\right)_{i^1_j,1} & \cdots & \prod_{j=1}^N \left(\bmA_{\sigma_j}\right)_{i^1_j,r} \\
                    \vdots & \vdots & \vdots \\
                    \prod_{j=1}^N \left(\bmA_{\sigma_j}\right)_{i^s_j,1} & \cdots & \prod_{j=1}^N \left(\bmA_{\sigma_j}\right)_{i^s_j,r}
                \end{array}
            \right]
        \right]_v
        \\
        &=
        \sum_{w=1}^s \left( \tilde{y}_{i_1^w,\dots,i_N^w} \left(\prod_{j=1}^N \left(\bmA_{\sigma_j}\right)_{i^w_j,v}\right) \right)
        \\
        &=
        \left[ \left(\odot_{j=1}^N \bmA_{\sigma_j}\right)^T \operatorname{vec}\left(\bctlY\right)\right]_v.
    \end{aligned}
\end{multline}
\qed

\bibliography{references}
\bibliographystyle{IEEEtran}

\vfill

\end{document}